\documentclass[12pt]{amsart}
\usepackage{somedefs,tracefnt,latexsym,,rawfonts,epsfig,amsfonts,amssymb,latexsym,enumerate,amscd}
\def\cal{\mathcal}
\def\Bbb{\mathbb}

\def\G{\Gamma}

\def\ul{\underline}
\def\F{\cal {FICWF}}
\def\flf{FICwF^{_{\underline L}{\cal H}^?_*}_{\cal {FIN}}}
\def\flv{FICwF^{_{\underline L}{\cal H}^?_*}_{\cal {VC}}}

\newtheorem{thm}{Theorem}[section]

\newtheorem{lemma}{Lemma}[section]
\newtheorem{cor}{Corollary}[section]

\newtheorem{defn}{Definition}[section]

\newtheorem{rem}{Remark}[section]
\numberwithin{equation}{section}
\begin{document}
\date{April 17, 2007.}
\title[$L$-theory of groups]{$L$-theory of groups with unstable derived series}
\author[S.K. Roushon]{S. K. Roushon}
\address{School of Mathematics\\
Tata Institute\\
Homi Bhabha Road\\
Mumbai 400005, India}
\email{roushon@math.tifr.res.in} 
\urladdr{http://www.math.tifr.res.in/\~\ roushon/}
\begin{abstract} In this short note we prove that the Farrell-Jones 
Fibered Isomorphism Conjecture in $L$-theory, after inverting $2$, is 
true for a group whose some derived subgroup is free.

\end{abstract} 



\maketitle


\section{Introduction} In \cite{R} it was 
shown that the 
Fibered Isomorphism Conjecture in 
$L^{-\infty}$-theory ([\cite{FJ}, 1.7]), 
after inverting $2$ (that is for 
$\ul L^{-\infty}=L^{-\infty}\otimes {\Bbb Z}[\frac{1}{2}]$),  
is true for a large class (denoted $\F$) 
of groups including poly-free 
groups and one-relator groups.

Here we deduce the conjecture, using the results in \cite{R},  
for groups whose some 
derived subgroup is free. See Remark \ref{remark} 
regarding the relevance of this class of groups.
 
Throughout the article, by `group' we mean 
`countable group' 

We prove the following.

\begin{thm} \label{solvable} Let $\G$ be a group. 
Then the Fibered Isomorphism Conjecture
of Farrell and Jones for the $\ul L^{-\infty}$-theory  
is true for the
group $\G\wr F$ if it is true 
for $\G^{(n)}\wr F$ for some $n$, where $F$ is a finite 
group and $\G^{(n)}$ denotes the $n$-th derived subgroup 
of $\G$.
 
In other words, the $\flv$ (or the $\flf$) is true 
for $\G$ if the $\flv$ (or the $\flf$) is true for 
$\G^{(n)}$ for some $n$.\end{thm}

For notations and statement of the conjecture 
see [\cite{R}, section 2].

\begin{proof}[Proof of Theorem \ref{solvable}] 
Consider the following exact sequence. 
$$1\to \G^{(n)}\to \G\to \G/\G^{(n)}\to 1.$$ 

Note that $\G/\G^{(n)}$ is a solvable group. Hence 
applying the hypothesis, 
Corollary \ref{solvcor} and [\cite{R}, 
$(2)$ of lemma 2.13] we complete the proof of the 
Theorem.\end{proof}

\begin{cor} \label{free} Let $G$ be a finite index subgroup of a 
group $\G$. 
Assume that $G^{(n)}$ is a free 
group for some $n$. Then the $\flv$ (or the $\flf$) is true 
for $\G$.\end{cor}

\begin{proof} Using Theorem \ref{mth} and Lemma 
\ref{imp} we can assume that $\G=G$. Next we only need to recall 
that by 
[\cite{R}, main lemma] the $\flv$ (or the $\flf$) 
is true for any free group and then apply Theorem 
\ref{solvable}.\end{proof}

\begin{rem}\label{remark}{\rm Groups whose derived series 
does not stabilize (or some derived subgroup 
is free or surjects onto a free group) 
are of interest in group theory and topology. 
See \cite{R1}. In fact in \cite{R1} we predicted that 
these kind of groups appear more often than  
other groups.}\end{rem}

Let us recall the following definition from \cite{R}.

\begin{defn} \label{begin} ([\cite{R}, definition 1.1]) 
{\rm Let $\F$ be the smallest
class of groups satisfying the following conditions

\begin{itemize}
\item The following groups belong to $\F$.

1. Finite groups. 2. Finitely generated free groups. 3. Cocompact
discrete subgroups of linear Lie groups with finitely many
components.
\item (Subgroup) If $H<G\in \F$ then $H\in \F$
\item (Free product) If $G_1, G_2\in \F$ then $G_1*G_2\in \F$.
\item (Direct limit) If $\{G_i\}_{i\in I}$ is a directed sequence
of groups with
$G_i\in \F$. Then the limit $\lim_{i\in I}G_i\in \F$.
\item (Extension) For an exact sequence of groups $1\to K\to G\to N\to 1$, if
$K, N\in \F$ then $G\in \F$.
\end{itemize}}\end{defn}

Let $A$ and $B$ be two groups then by definition 
the wreath product $A\wr B$ is the semidirect product $A^B\rtimes B$ 
where the action of $B$ on $A^B$ is the regular action. 
Let ${\cal {VC}}$ and 
${\cal {FIN}}$ denote the class of virtually cyclic groups 
and the class
of finite groups respectively.

We proved the following theorem in [\cite{R}, theorem 1.1].

\begin{thm} \label{mth} ([\cite{R}, theorem 1.1]) 
Let $\Gamma\in \F$. 
Then the following assembly maps
are isomorphisms for all $n$, for any group homomorphism $\phi :G\to
\Gamma\wr F$ and for any finite group $F$. 
$${\cal H}^G_n(p, {\bf \ul L}^{-\infty}):{\cal
H}^G_n(E_{\phi^*{\cal
{VC}}(\G\wr F)}(G), {\bf \ul L}^{-\infty})\to
{\cal H}^G_n(pt, {\bf \ul L}^{-\infty})\simeq
\ul L_n^{-\infty}({\Bbb
Z}G).$$
$${\cal H}^G_n(p, {\bf \ul L}^{-\infty}):{\cal
H}^G_n(E_{\phi^*{\cal
{FIN}}(\G\wr F)}(G), {\bf \ul L}^{-\infty})\to
{\cal H}^G_n(pt, {\bf \ul L}^{-\infty})\simeq
\ul L_n^{-\infty}({\Bbb
Z}G).$$

In other words the Fibered Isomorphism Conjecture
of Farrell and Jones for the $\ul L^{-\infty}$-theory  
is true for the
group $\G\wr F$. Equivalently, 
the $\flv(\G)$ and the $\flf(\G)$ are satisfied 
(see [\cite{R}, definition 2.1] for notations).\end{thm} 

In \cite{R} we showed that $\F$ contains some well-known 
classes of groups. Here we see that $\F$ also contains 
any virtually solvable group. 

\begin{thm} \label{hyper} $\F$ contains the class of  
virtually solvable groups.\end{thm}   

\begin{proof} Let $\G$ be a virtually solvable group. 
Using the `direct limit' condition in the definition 
of $\F$ we can assume that 
$\G$ is finitely generated, for any countable infinitely 
generated group is a direct limit of finitely generated 
subgroups.
The following Lemma shows that we can also assume that the 
group $\G$ is solvable.

\begin{lemma} \label{imp} 
Let $G$ be a finitely generated group and   
contains a finite index subgroup $K$. If $K\in \F$ 
then $G\in \F$.\end{lemma}

\begin{proof} By taking the intersection of all 
conjugates of $K$ in $G$ we get a subgroup $K'$ 
of $G$ which is normal and of finite index in 
$G$. Therefore, we can use `subgroup' and `extension' 
conditions in the definition of $\F$ to conclude 
the proof of the Lemma.\end{proof}

Hence we have $\G$ a finitely generated solvable 
group. We say that $\G$ is $n$-step 
solvable if $\G^{(n+1)}=(1)$ and $\G^{(n)}\neq (1)$. The 
proof is by induction on $n$. Since 
countable abelian groups belong to $\F$ 
(see [\cite{R}, lemma 4.1]), the 
induction starts. 

So assume that a finitely generated $k$-step solvable 
group for $k\leq n-1$ belong to $\F$ and $\G$ is 
$n$-step solvable.

We have the following exact 
sequence. $$1\to \G^{(n)}\to \G\to \G/\G^{(n)}\to 1.$$

Note that $\G^{(n)}$ is abelian and $\G/\G^{(n)}$ is 
$(n-1)$-step solvable. Using the `extension' condition 
and the induction hypothesis we complete the proof.
\end{proof}

Applying Theorems \ref{mth} and \ref{hyper} we get 
the following.

\begin{cor} \label{solvcor} The $\flv$ (or the $\flf$) 
is true for any virtually solvable group.\end{cor} 

\newpage
\bibliographystyle{plain}

\begin{thebibliography}{60}

\bibitem{FJ}
F.T. Farrell and L.E. Jones,
\newblock Isomorphism conjectures in algebraic $K$-theory, 
\newblock {\em J. Amer. Math. Soc.}, 6 (1993), 249-297.

\bibitem{LR}
W. L\"{u}ck and H. Reich,
\newblock The Baum-Connes and the Farrell-Jones Conjectures in $K$- and
$L$-theory,
\newblock  Handbook of K-theory Volume 2, editors: E.M. Friedlander, D.R.
Grayson, (2005) 703-842, Springer.

\bibitem{R} 
S.K. Roushon,
\newblock The isomorphism conjecture in $L$-theory: poly-free 
groups and one-relator groups, math.KT/0703879.

\bibitem{R1}
\bysame,
\newblock Topology of 3-manifolds and a class of groups II,
\newblock {\em Bol. Soc. Mat. Mexicana (3)} 10 (2004), Special Issue,
467-485.

\end{thebibliography}
\ifx\undefined\bysame
\newcommand{\bysame}{\leavevmode\hbox to3em{\hrulefill}\,}
\fi

\medskip

\end{document}